\documentclass[11pt]{amsart}
\usepackage{graphics,latexsym,amssymb,amsmath}
\usepackage[dvips]{graphicx}
\usepackage{psfrag}
\input{amssym.def} 
\graphicspath{{figs/}} 

\newtheorem{thm}{Theorem}[section]
\newtheorem{dfn}[thm]{Definition}

\newtheorem{rmk}[thm]{Remark}
\newtheorem{rmks}[thm]{Remarks}
\newtheorem{cor}[thm]{Corollary}
\newtheorem{prop}[thm]{Proposition}
\newtheorem{lem}[thm]{Lemma}

\newcommand{\EPf}{\hbox{}\hfill$\Box$\vspace{.5cm}}

\newcommand{\R}{{{\mathbb R}}}

\newcommand{\F}{{{\mathcal F}}}

\newcommand{\E}{{{\goth E}}}

\newcommand{\Id}{{\rm Id}}
\newcommand{\Aut}{{\rm Aut}}
\newcommand{\Cay}{{\rm Cay}}
\newcommand{\spann}{{\rm span}}
\newcommand{\diag}{{\rm diag}}
\newcommand{\supp}{{\rm supp}}

\date{}
\begin{document}
\title[Spectral representations of vertex transitive graphs]{Spectral
  representations of vertex transitive graphs, Archimedean solids and
  finite Coxeter groups} 
\author{Ioannis Ivrissimtzis and Norbert
  Peyerimhoff}
\address{Durham University, DH1 3LE, Great Britain}
\email{\{ioannis.ivrissimtzis,norbert.peyerimhoff\}@durham.ac.uk}
\date{}
\subjclass[2000]{20F55 (Primary) 05C10, 05C50, 05C62, 05C81, 52B15 (Secondary)}

\begin{abstract}
  In this article, we study eigenvalue functions of varying transition
  probability matrices on finite, vertex transitive graphs. We prove
  that the eigenvalue function of an eigenvalue of fixed higher
  multiplicity has a critical point if and only if the corresponding
  spectral representation is equilateral. We also show how the
  geometric realisation of a finite Coxeter group as a reflection
  group can be used to obtain an explicit orthogonal system of
  eigenfunctions. Combining both results, we describe the behaviour of
  the spectral representations of the second highest eigenvalue
  function under the change of the transition probabilities in the
  case of Archimedean solids.
\end{abstract}

\maketitle

\section{Introduction and statement of results}

The main objects of interest in this paper are {\em spectral
representations} associated to random walks on finite graphs (see
Sections \ref{subsec:basics} and \ref{subsec:specrep} for the
definitions). We consider the particular case of vertex transitivity,
which comprises the large class of Cayley graphs. In our main result
(Theorem \ref{thm:main} below), we prove a correspondence between the
critical points of an eigenvalue function (under the change of the
invariant transition probabilities) and the points where the
associated spectral representation is equilateral. In Sections
\ref{subsec:coxres} and \ref{subsec:archsol}, we specialise our
considerations to finite Coxeter groups and one-skeleta of Archimedean
solids.

\subsection{Basic graph theoretical notation} \label{subsec:basics}
Let $G = (V,E)$ be a finite, simple (i.e., no loops and multiple
edges) graph with vertex set $V = \{1,\dots,n\}$ and set of undirected
edges $E$. An edge is represented by a set $\{i,j\} \subset V$ with $i
\neq j$. A (time reversible) random walk on $G$ is given by a {\em
  symmetric} stochastic matrix $P = (p_{ij}) \in \R^{n \times n}$,
where $p_{ij}$ is the transition probability from vertex $i$ to vertex
$j$. For $i \neq j$, we require $p_{ij} = 0$ if $\{i,j\} \not\in E$.
Even though there are no loops in $G$, we allow the diagonal elements
$p_{ii}$ to be positive. ($p_{ii}$ represents the probability for the
random walk to stay at the vertex $i$.) The set of all matrices $P$ of
the above type are a convex subset of $\R^{n \times n}$, which we
denote by $\Pi_G$. We think of a matrix $P \in \Pi_G$ as a linear
operator on the vector space $l^2(G)$ of (real-valued) functions on
the vertices, i.e.,
$$ P f(i) = p_{ii} f(i) + \sum_{j \sim i} p_{ij} f(j), $$
where $j \sim i$ means that $\{i,j\} \in E$. The inner product on
$l^2(G)$ is given by $\langle f,g \rangle = \sum_{i=1}^n f(i)
g(i)$. Let $\sigma(P)$ denote the spectrum of $P$ with eigenvalues
$$ 1=\lambda_0(P) \ge \lambda_1(P) \ge \cdots \lambda_{n-1}(P) \ge -1, $$
counted with multiplicity.  Let $f_0(i) = \frac{1}{\sqrt{n}}$. The
Rayleigh quotient representation of the second highest eigenvalue function
\begin{equation} \label{eq:rayleigh} 
  \lambda_1(P) = \sup_{f \bot f_0 \mid f \neq 0 } \frac{\langle P f,f \rangle}
  {\Vert f \Vert^2}
\end{equation}
implies that $\lambda_1: \Pi_G \to [-1,1]$ is {\em convex} (see the
proof of Proposition \ref{prop:lambda1boundary} in Section
\ref{sec:genresults}). The functions $\lambda_i: \Pi_G \to [-1,1]$ are
continuous (see, e.g., \cite[Theorem (1,4)]{Marden2005}), but these
functions fail to be analytic at those points where eigenvalues of
higher multiplicity bifurcate. We refer the reader to, e.g.,
\cite[Chapter 2]{Kato95}, for more information about these subtle
regularity issues. The special operator $P = (p_{ij})$ with vanishing
main diagonal ($p_{ii} = 0$ for all $i \in V$), and for which all
other transition probabilities $p_{ij}$ are equal to $1/{\rm deg}(i)$,
is called the {\em canonical Laplacian}.

\subsection{Spectral representations} \label{subsec:specrep}

The idea of a {\em spectral representation} is to use a higher
multiplicity eigenvalue of the matrix $P$ to obtain a ''geometric
realisation'' of the combinatorial graph $G$ in Euclidean
space. Assume that $\lambda \in \sigma(P)$ is an eigenvalue of $P$ of
multiplicity $k$, and $\phi_1,\dots,\phi_k$ is an orthonormal base of
eigenfunctions of the eigenspace $\E_\lambda(P)$. The corresponding
spectral representation is the map
$$ \Phi = \Phi_{P,\lambda}: V \to \R^k, \quad \Phi(i) = 
(\phi_1(i),\dots,\phi_k(i)), $$ 
i.e., the simultaneous evaluation of all eigenfunctions at a given
vertex. The spectral representation depends on the choice of the
orthonormal base only up to an orthonormal transformation in $\R^k$.

There are often striking geometric and spectral analogies between
the discrete setting of graphs and the smooth setting of Riemannian
manifolds. In the context of Riemannian manifolds, the simultaneous
evaluation of eigenfunctions of the Laplacian were considered, for
example, in the so-called {\em nice (minimal isometric) embeddings} of
strongly harmonic manifolds into Euclidean spheres (see \cite[Chapter
6G]{Besse78}).

\begin{dfn}
  A spectral representation $\Phi: V \to \R^k$ is {\em faithful} if
  $\Phi$ is injective. It is {\em equilateral} if all images of edges
  have the same Euclidean length, i.e.,
  $$ \Vert \Phi(i_1)-\Phi(j_1) \Vert = \Vert \Phi(i_2)-\Phi(j_2) \Vert, $$
  for all pairs of edges $\{i_1,j_1\}, \{i_2,j_2\} \in E$, where $\Vert
  \cdot \Vert$ denotes the Euclidean norm.
\end{dfn}

A particularly strong faithfulness result for $3$-connected, planar
graphs in the case that the second highest eigenvalue has multiplicity
three was obtained in \cite{LSch-99}.

\subsection{Vertex transitive graphs}

In this paper, we focus on finite {\em vertex transitive graphs},
i.e., we assume that the automorphism group ${\rm Aut}(G)$ acts
transitively on the vertex set $V$. Particular examples of vertex
transitive graphs are Cayley graphs of groups. Below, we introduce
equivalence classes of edges and, to have enough flexibility, we
consider subgroups $\Gamma \subset {\rm Aut}(G)$ which still act
transitively on the vertices. We define a $\Gamma$-action on the
space $\Pi_G$ of matrices as follows:
$$ (\gamma P)_{ij} = p_{\gamma i,\gamma j} \qquad 
\text{for all $P = (p_{ij}) \in \Pi_G$.} $$ 
A random walk and its corresponding matrix $P \in \Pi_G$ is called {\em
  $\Gamma$-invariant}, if $\gamma P = P$ for all $\gamma \in \Gamma$.
Note that the main diagonal $(p_{11},\dots,p_{nn})$ of every
$\Gamma$-invariant matrix $P$ is constant. The large automorphism group of
a vertex transitive graph makes the occurence of eigenvalues of higher
multiplicities for $\Gamma$-invariant matrices more likely, and it is
natural to make use of connections between these eigenvalues and the
representation theory of $\Gamma$.

The group $\Gamma$ induces an equivalence relation on the set of
edges: $\{ i,j \} \in E$ is equivalent to all edges $\{ \gamma i,
\gamma j \}$ with $\gamma \in \Gamma$. The {\em multiplicity} of an
equivalence class $[e] \subset E$ is the number of edges in $[e]$
meeting at the same vertex. Let $[e_1],\dots,[e_N]$ be the $\Gamma$
equivalence classes of edges and $m_1,\dots,m_N$ be its
multiplicities. The set of {\em $\Gamma$-invariant matrices in $\Pi_G$
  with vanishing main diagonal} is a convex subset, which we identify
with the simplex
\begin{equation} \label{eq:deltagen}
\Delta_\Gamma = 
\{ (x_1,\dots,x_N) \in [0,1]^N \mid \sum_j m_j x_j = 1 \}.
\end{equation}
The point $X = (x_1,\dots,x_N) \in \Delta_\Gamma$ corresponds to the
matrix $P_X = (p_{ij})$, given by
$$ p_{ij} = \begin{cases} 0, & \text{if $i=j$ or $\{i,j\} \not\in E$,} \\
  x_k, & \text{if $\{i,j\} \in [e_k]$.} \end{cases} $$ 
For $P = (p_{ij}) \in \Pi_G$, let $G_P = (V,E_P)$ denote the
subgraph of $G$ with edges $E_P = \{ \{i,j\} \in E \mid p_{ij} > 0
\}$. Then, for every interior point $X \in {\rm int}(\Delta_\Gamma)$, we
have $G_{P_X} = G$ (since the entries $p_{ij}$ associated to all edges
$\{i,j\}$ are strictly positive), and the spectrum $\sigma(P_X)$
is symmetric with respect to the origin if and only if $G_{P_X}$ is
bipartite.

Let us now discuss the special case of Cayley graphs. A finite
symmetric set $S \subset \Gamma$ of generators of a group $\Gamma$ is
called {\em minimal} if for every $s \in S$, $S - \{s,s^{-1}\}$ is no
longer a set of generators. The Cayley graph of $\Gamma$ with respect
to $S$ is denoted by $\Cay(\Gamma,S)$, its vertices are the group
elements, i.e., $V = \Gamma$, and two vertices $\gamma, \gamma'$ are
connected by an edge if and only if $\gamma' = \gamma s$ for some $s
\in S$. If $S = \{s_1,\dots,s_r\}$ is a minimal symmetric set of
generators, we distinguish the generators of order $2$ (since they
appear only once in $S$) from the ones with higher order, by rewriting
them as
\begin{equation} \label{eq:genset}
S = \{s_1,\dots,s_\nu,\tau_1^{\pm 1},\dots,\tau_\mu^{\pm 1} \}, 
\end{equation}
with $\nu + 2 \mu = r$. Note that the edges $\{e,\tau_j\}$ and
$\{e,\tau_j^{-1}\}$ are equivalent, and the corresponding simplex
$\Delta_\Gamma$ is given by
\begin{equation} \label{eq:gensimplex}
\Delta_\Gamma = \{(x_1,\dots,x_{\nu+\mu}) \in [0,1]^N \mid \sum_{j=1}^\nu x_j 
+ 2 \sum_{j=1}^\mu x_{\nu+j} = 1\}. 
\end{equation}
The following facts follow from the convexity of $\lambda_1: \Pi_G \to
[-1,1]$ (see Section \ref{sec:genresults} for the proof).

\begin{prop} \label{prop:lambda1boundary} 
  Let $G = (V,E)$ be a finite, connected, simple graph and $\Gamma
  \subset \Aut(G)$ be vertex transitive. Then a global minimum of
  $\lambda_1: \Pi_G \to [-1,1]$ is assumed at a matrix in
  $\Delta_\Gamma$. 

  If $G = \Cay(\Gamma,S)$ is the Cayley graph of a finite group
  $\Gamma$ with respect to a minimal symmetric set $S$ of generators,
  then we have
  \begin{equation} \label{eq:convbound} 
  \lim_{n \to \infty} \lambda_1(P_{X_n}) = 1 
  \end{equation}
  for every sequence $X_n \to \partial \Delta_\Gamma$, and a global
  minimum of $\lambda_1$ is assumed at an interior point of $\Delta_\Gamma$.
\end{prop}

Note that the above result does not rule out that $\lambda_1$ may also
have other global minima at matrices $P \in \Pi_G - \Delta_\Gamma$.

Our main general result is the following relation between critical
points of eigenvalue functions and equilateral spectral
representations:

\begin{thm} \label{thm:main} 
  Let $G = (V,E)$ be a finite, connected, simple graph and $\Gamma
  \subset \Aut(G)$ be vertex transitive. Let $U \subset \Delta_\Gamma$
  be an open set and $\lambda: U \to [-1,1]$ be a smooth function such
  that $\lambda(X):=\lambda(P_X)$ is an eigenvalue of $P_X$ with fixed
  multiplicity $k \ge 2$ for all $X \in U$. $X_0 \in U$ is a critical
  point of the function $\lambda$ if and only if the spectral
  representation $\Phi = \Phi_{P_{X_0},\lambda(X_0)}: V \to S^{k-1}$
  is equilateral.
\end{thm}

It is shown in Lemma \ref{lem:basic} (see Section
\ref{sec:genresults}) that, for vertex transitive graphs, the image of
every $\Gamma$-invariant spectral representation
$\Phi=\Phi_{P_X,\lambda}$ (with $\lambda \in \sigma(P_X)$ of
multiplicity $k$) lies on an Euclidean sphere $S^{k-1} \subset \R^k$, and
that equivalent edges are mapped to segments with the same Euclidean
length, i.e., $\Vert \Phi(\gamma i) - \Phi(\gamma j) \Vert = \Vert
\Phi(i) - \Phi(j) \Vert$. The above theorem states that at critical
points of the eigenvalue function all Euclidean images of edges have
the same length (not only the equivalent ones). 

\begin{rmks}
(a) Special examples of critical points are minima of a smooth
function. As another example for similarities between graphs and
Riemannian manifolds, we like to mention the following result in
Riemannian geometry: The first non-zero Laplace-eigenvalue
$\lambda_1(M) > 0$ of a closed Riemannian manifold $(M,g)$ of
dimension $n$ with lower positive Ricci curvature is minimal if and
only if $M$ is isometric to the $n$-dimensional round sphere (Obata's
theorem, see \cite{Ob62} or \cite{BGM71}). Here we also
have the phenomenon that a critical point of the eigenvalue function
is assumed in the case of a very symmetric geometry.

(b) For extremal eigenvalues of the Laplace matrix of general graphs,
related embedding interpretations arose, e.g., in \cite{GHW08,GHR10}
in studying the semidefinite duals of associated eigenvalue
optimization problems. The relation of these results to the vertex
symmetric graphs studied here becomes more apparent when symmetry is
exploited in the corresponding optimization problems by the techniques
described, e.g., in \cite{deK10,BGSV10}. The precise nature of this
relation, however, still needs to be explored further.
\end{rmks}

Standard arguments in representation theory yield the following useful result:

\begin{prop} \label{prop:main2} 
  Let $\Gamma$ be a finite group with a minimal symmetric set of
  generators $S$ given by \eqref{eq:genset} and $G = \Cay(\Gamma,S)$
  be the associated Cayley graph with the corresponding simplex
  $\Delta_\Gamma$ as in \eqref{eq:gensimplex}.
 
  Let $\rho: \Gamma \to O(k)$ be an irreducible representation,
  $\pi_r: \R^k \to \R$ be the projection to the $r$-th coordinate and
  $S^{k-1} \subset \R^k$ be the unit sphere. Let $p \in S^{k-1}$,
  $\lambda \in \R$, and $X = (x_1,\dots,x_{\nu+\mu}) \in \Delta_\Gamma$ such
  that
  \begin{equation} \label{eq:lamP} \lambda p = \sum_{j=1}^\nu x_j
    \rho(s_j)p + \sum_{j=1}^\mu x_{\nu+j} (\rho(\tau_j)p +
    \rho(\tau_j^{-1})p).
  \end{equation}
  Then the functions
  $$ \phi_r: \Gamma \to \R, \quad \phi_r(\gamma) := \pi_r(\rho(\gamma) p), 
  \quad 1 \le r \le k $$ 
  are pairwise orthogonal eigenfunctions of $P_X$ for the eigenvalue
  $\lambda$ satisfying $\Vert \phi_r \Vert^2 = \frac{|\Gamma|}{k}$.
\end{prop}

\begin{rmks} (a) This result implies that if the eigenspace
  $\E_\lambda(P_X)$ is an irreducible representation of $\Gamma$
  (i.e., $\phi_1,\dots,\phi_k$ span the whole eigenspace
  $\E_\lambda(P_X)$), then the associated spectral representation
  $\Phi: \Gamma \to S^{k-1}$ coincides with the orbit map
  $\Phi(\gamma) = \rho(\gamma) p_0$ of the rescaled point $p_0 =
  \sqrt{\frac{k}{|\Gamma|}}\, p \in \R^k$. Thus a natural question
  is whether eigenspace representations are irreducible, or whether
  different representations appear with the same eigenvalue. 
  
  (b) It can be shown, in the weaker case of a non-orthogonal
  irreducible representation $\rho: \Gamma \to GL(k,\R)$, that the
  functions $\phi_r$ are still a family of {\em linear independent}
  eigenfunctions of $P_X$.
\end{rmks}

\subsection{Finite irreducible Coxeter groups} \label{subsec:coxres}

Let us now consider the special case of a finite irreducible Coxeter
group $\Gamma = \langle S=\{s_1,\dots,s_k\} \mid
(s_is_j)^{m_{ij}}=e\rangle$ of rank ${\rm rk}(\Gamma)=k$ with
$m_{ij}=m_{ji} \ge 2$ and $m_{ii}=2$, i.e., $s_i^2 =e$. It was
suggested in \cite[Problem 10.8.7]{Lub2010} to study the eigenvalues
(or at least $\lambda_1$) of the canonical Laplacian for Coxeter
groups. Bacher \cite{Ba94} identified $\lambda_1$ of the canonical
Laplacian for symmetric groups. For the canonical Laplacian on
arbitrary finite Coxeter groups, Akhiezer \cite{Akh07} found an
explicit set of eigenvalues and a lower bound on their multiplicity in
case of irreducibility. The spectral gap of the canonical Laplacian
and the Kazdhan constant of all finite Coxeter groups was explicitly
derived by Kassabov in \cite[Section 6.1]{Kas09}. For infinite Coxeter
groups, it was proved in \cite{BJS88} that they do not have Kazdhan
property $(T)$. In this section we are concerned with Laplacians on
finite, irreducible Coxeter groups with variable weights.

Let $\Gamma \hookrightarrow O(k)$ be the geometric realisation of
$\Gamma$ as finite reflection group. The associated Cayley graph $G =
\Cay(\Gamma,S)$ is bipartite, since all relations of a Coxeter group
have even length. Let $\sigma_j \in O(k)$ be the reflections
corresponding to the generators $s_j$ and $n_1,\dots,n_k$ be the
associated simple roots. Let
\begin{equation} \label{eq:Pj}
p_j = (-1)^{j-1} n_1 \times \cdots \times \widehat{n_j} \times \cdots
\times n_k, 
\end{equation}
where $v_1\times \cdots \times v_{k-1}$ denotes the $(k-1)$-ary
analogue of the cross product in $\R^k$, and the hat over $n_j$ in
\eqref{eq:Pj} means that this term is dropped.  Then the open cone
$$ \F := \{ \alpha_1 p_1 + \cdots + \alpha_k p_k \mid \alpha_1,\dots,\alpha_k 
> 0 \} \subset \R^k, $$ is a fundamental domain of the $\Gamma$-action
on $\R^k$. $\Gamma$ preserves the unit sphere $S^{k-1}$, and a
spherical fundamental domain is given by $\F_0 = \F \cap S^{k-1}$. Let
$V = \det(n_1,\dots,n_k)$. Without loss of generality, we can assume
that $V > 0$, for otherwise we simply permute the set of
generators. The following result is a consequence of Proposition
\ref{prop:main2}.

\begin{cor} \label{cor:coxgroup} Let $\Gamma$ be a finite, irreducible
  Coxeter group and $\F_0 \subset S^{k-1}$ and $\Delta=\Delta_\Gamma$
  be as above. Then there exists smooth maps $\Psi_\Delta: \F_0 \to
  {\rm int}(\Delta)$ and $\Psi_\lambda: \F_0 \to (0,1)$, with
  $\Psi_\Delta$ bijective, such that, for every $p = \sum \alpha_j p_j
  \in \F_0$, the functions $\phi_i(\gamma) := \pi_i(\gamma p)$ are
  pairwise orthogonal eigenfunctions of $P_X$ on $\Cay(\Gamma,S)$ for
  the eigenvalue $\lambda=\Psi_\lambda(p)$, where $X=\Psi_\Delta(p)$.
  Moreover, $\Vert \phi_i \Vert^2 = \frac{|\Gamma|}{k}$ and the
  composition $\Psi_\lambda \circ \Psi_\Delta^{-1}: {\rm int}(\Delta)
  \to (0,1)$ is analytic. The simultaneous evaluation
  $$ \Phi(\gamma): \Gamma \to S^{k-1}, \quad \Phi(\gamma) = 
  (\phi_1(\gamma),\dots,\phi_k(\gamma)) = \gamma p $$ 
  is faithful, and the Euclidean lengths of the images of equivalence
  classes of edges under $\Phi$ are given by
  $$ \Vert p - \sigma_j(p) \Vert = 2 \alpha_j V. $$
\end{cor}

\begin{rmk}
  The explicit description of the maps $\Psi_\Delta$, $\Psi_\lambda$
  and the composition $\Psi_\lambda \circ \Psi_\Delta^{-1}$ is given
  by the equations \eqref{eq:calcX'}, \eqref{eq:calcl'},
  \eqref{eq:PsisP} and \eqref{eq:PsiPsi} in Section
  \ref{sec:coxresults1}.
\end{rmk}

The next result follows from a slight modification of a calculation given
in Kassabov \cite{Kas09}.

\begin{prop} \label{prop:kassabov} Let $\Gamma, \Psi_\Delta$ and
  $\Psi_\lambda$ be as in Corollary \ref{cor:coxgroup}. Then the map
  $\Psi_\lambda \circ \Psi_\Delta^{-1}: {\rm int}(\Delta) \to (0,1)$
  coincides with the second highest eigenvalue function $\lambda_1:
  {\rm int}(\Delta) \to (0,1)$. Consequently, the second highest eigenvalue
  $\lambda_1(X)$ of $P_X$ has always multiplicity $\ge {\rm rk}(\Gamma)$. 
\end{prop} 

The proof of the next result on the exact multiplicity of the second
highest eigenvalue for particular Coxeter groups is based on elegant
arguments of van der Holst \cite{vdH95}. He used these arguments to
give a direct combinatorial proof of Colin de Verdier{\`e}'s planarity
characterisation ''$\mu(G) \le 3$''.
 
\begin{prop} \label{prop:A3B3H3} Let $\Gamma$ be one of the Coxeter
  groups $A_3$, $B_3$ or $H_3$.  Then the second highest eigenvalue
  $\lambda_1(X)$ of $P_X$ has multiplicity equals three for all $X \in
  {\rm int}(\Delta)$.
\end{prop}

\begin{rmk}
  (a) The heart of the proof of Proposition \ref{prop:A3B3H3}, namely
  van der Holst's argument, is geometric and depends on the planarity
  of the associated Cayley graphs. It is likely that for every finite,
  irreducible Coxeter group $\Gamma$ (not only $A_3, B_3, H_3$) the
  multiplicity of the second highest eigenvalue function is constant
  and equal to the rank of $\Gamma$. The techniques in Kassabov's
  paper \cite{Kas09} might be useful to prove this general statement.

  (b) The value of $\lambda_1(X)$ has a well known {\em dynamical
    interpretation}: Our Cayley graphs are bipartite, i.e., we have a
  partition $V = V_0 \cup V_1$. $\lambda_1(X)$ measures the
  convergence rate of the corresponding random walk to the
  equidistribution {\em (mixing rate)} on each set of vertices $V_i$
  under even time steps (even time steps are needed because of the
  bipartiteness). The validity of the multiplicity assumption in (a)
  together with our main result (Theorem \ref{thm:main}) would allow
  us to explicitly determine, for all finite, irreducible Coxeter
  groups, the transition probabilities of a random walk with the {\em
    fastest mixing rate} on the corresponding Cayley graphs. In fact,
  this is precisely how we will prove Theorem \ref{thm:seceigvalarch}
  below.
\end{rmk}

\subsection{Archimedean solids} \label{subsec:archsol}

The Cayley graph of the Coxeter groups $A_3, B_3$ and $H_3$ (with
respect to their set of standard generators $\{s_1, s_2, s_3\}$)
conincide, combinatorially, with the one-skeleta of the Archimedean
solids with the vertex configurations $(4,6,6), (4,6,8)$, and
$(4,6,10)$, respectively.

{\em Archi\-me\-dean solids} are polyhedra in $\R^3$ such that all
faces are regular polygons, and which have a symmetry group acting
transitively on the vertices. (Note, however, that the prisms,
antiprisms and Platonic solids, which also have these properties, are
excluded). The $13$ Archimedean solids are classified via their {\em
  vertex configurations}: The vertex configuration $(m,n,k)$ stands
for the solid where an $m$-gon, an $n$-gon and a $k$-gon (in this
order) meet at every vertex. We will use this notation also for
Platonic solids (e.g., the icosahedron is denoted by
$(3,3,3,3,3)$). The spectra of the canonical Laplacians (on the
one-skeleta) of all Archimedean solids were explicitly calculated in
\cite{ST-98}. For all these graphs, the second highest eigenvalue of
the canonical Laplacian has multiplicity three. The corresponding
spectral representation is faithful and represents a polyhedron in
$\R^3$ (this follows, e.g., from the general result in
\cite{LSch-99}), but this polyhedron is generally {\em not
  equilateral}. It is natural to study the deformation of this
polyhedron under changes of the $\Gamma$-invariant transition
probabilities (assuming that the multiplicity of $\lambda_1$ does not
change), and to find points at which the spectral representation is
equilateral.

We will carry this out in the case of the largest Archimedean solid,
namely the truncated icosidodecahedron $(4,6,10)$. We will also
explain, how the corresponding results read in the case of the
Archimedean solids $(4,6,8)$ and $(4,6,6)$. The proofs for these cases
are completely analogous.

Let $G=(V,E)$ be the one-skeleton of the Archimedean solid $(4,6,10)$.
The automorphism group of $G$ is the full icosahedral group and acts
simply transitively on the vertex set $V$, and is isomorphic to
$H_3$. Considering $G$ as a planar graph, its faces are $4$-, $6$- and
$10$-gons. $G$ is $3$-connected, has $120$ vertices and every vertex
has degree three (see Figure \ref{fig:specrep111} below). Let $\varphi =
\frac{1+\sqrt{5}}{2}=2 \cos \frac{\pi}{5}$ be the golden ratio. Our
previous results imply the following facts for $\lambda_1$.

\begin{thm} \label{thm:seceigvalarch} 
  Let $G$ be the $1$-skeleton of the Archimedean solid $(4,6,10)$ and
  $\Gamma = \Aut(G)$. The simplex of $\Gamma$-invariant transition
  probabilities is
  $$ \Delta = \Delta_\Gamma = \{ (x,y,z) \mid x,y,z \ge 0 , x+y+z=1 \}, $$
  where $x,y,z$ are the transition probabilities for the
  edge-equivalence classes separating $4$- and $6$-gons, $4$- and
  $10$-gons, and $6$- and $10$-gons, respectively. Then the
  restriction of $\lambda_1: \Pi_G \to [-1,1]$ to ${\rm int}(\Delta)
  \subset \Pi_G$ is analytic and strictly convex, and $\lambda_1(X)$
  has multiplicity three for all $X \in {\rm int}(\Delta)$.  Moreover,
  $X_0 = \frac{1}{14+5\varphi}(5,3+3\varphi,6+2\varphi)$ is the unique
  point in $\Delta$ at which $\lambda_1$ assumes its global minimum
  with
  $$ \lambda_1(X_0) = \frac{10+7\varphi}{14+5\varphi}. $$
  The corresponding spectral representation $\Phi_{X_0}: V \to S^2$ is
  faithful and equilateral.
\end{thm}

Let us stress, again, that for $X_* = (1/3,1/3,1/3) \in
\Delta$, the above Theorem implies that the spectral representation of
$P_{X_*}$ for $\lambda_1(X_*)$ does {\em not} reproduce
the Archimedean solid, one has to choose the point $X_0 \in \Delta$
instead (see Figure \ref{fig:specrep111}).

\begin{figure}[h]
  \begin{center}
  \makebox[6cm]{\includegraphics[width=6cm]{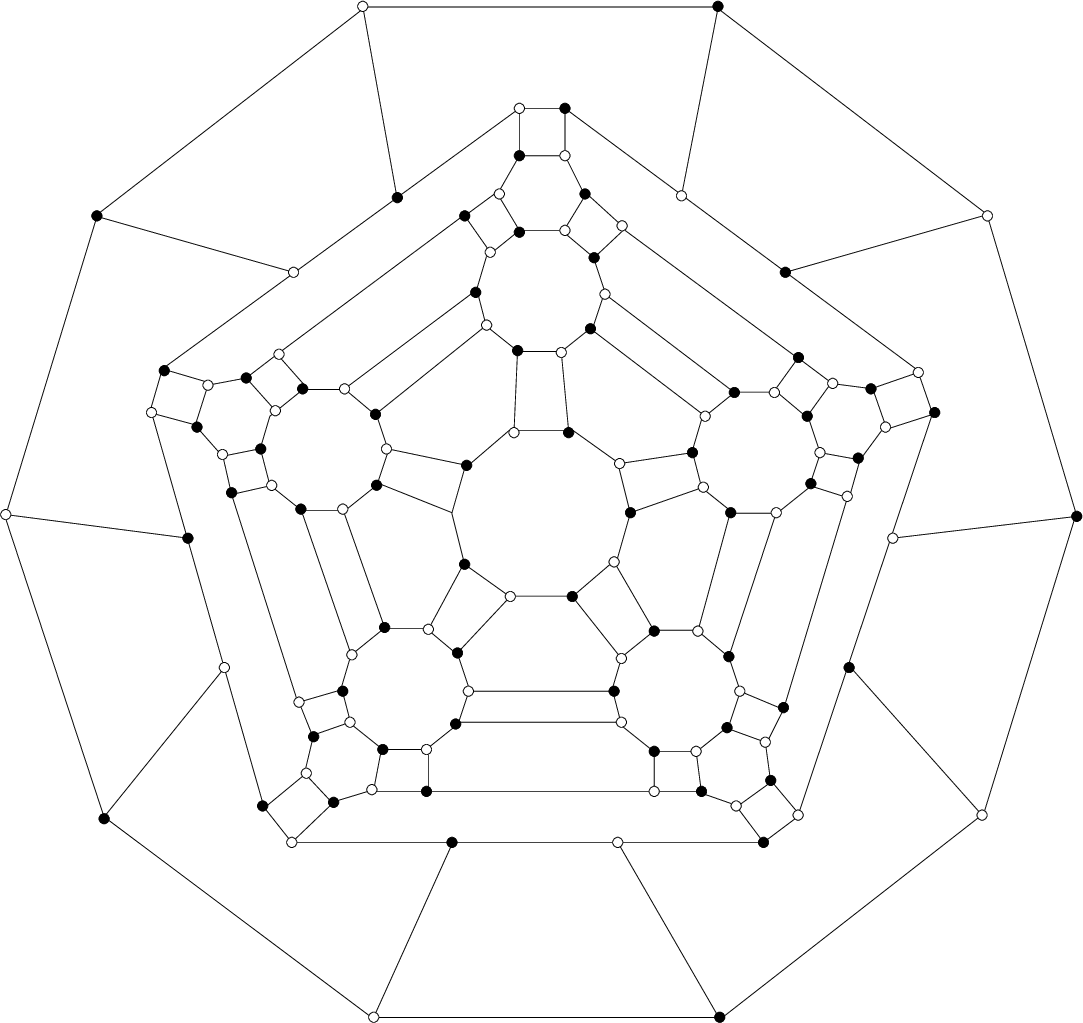}}
  \makebox[6cm]{\includegraphics[width=6cm]{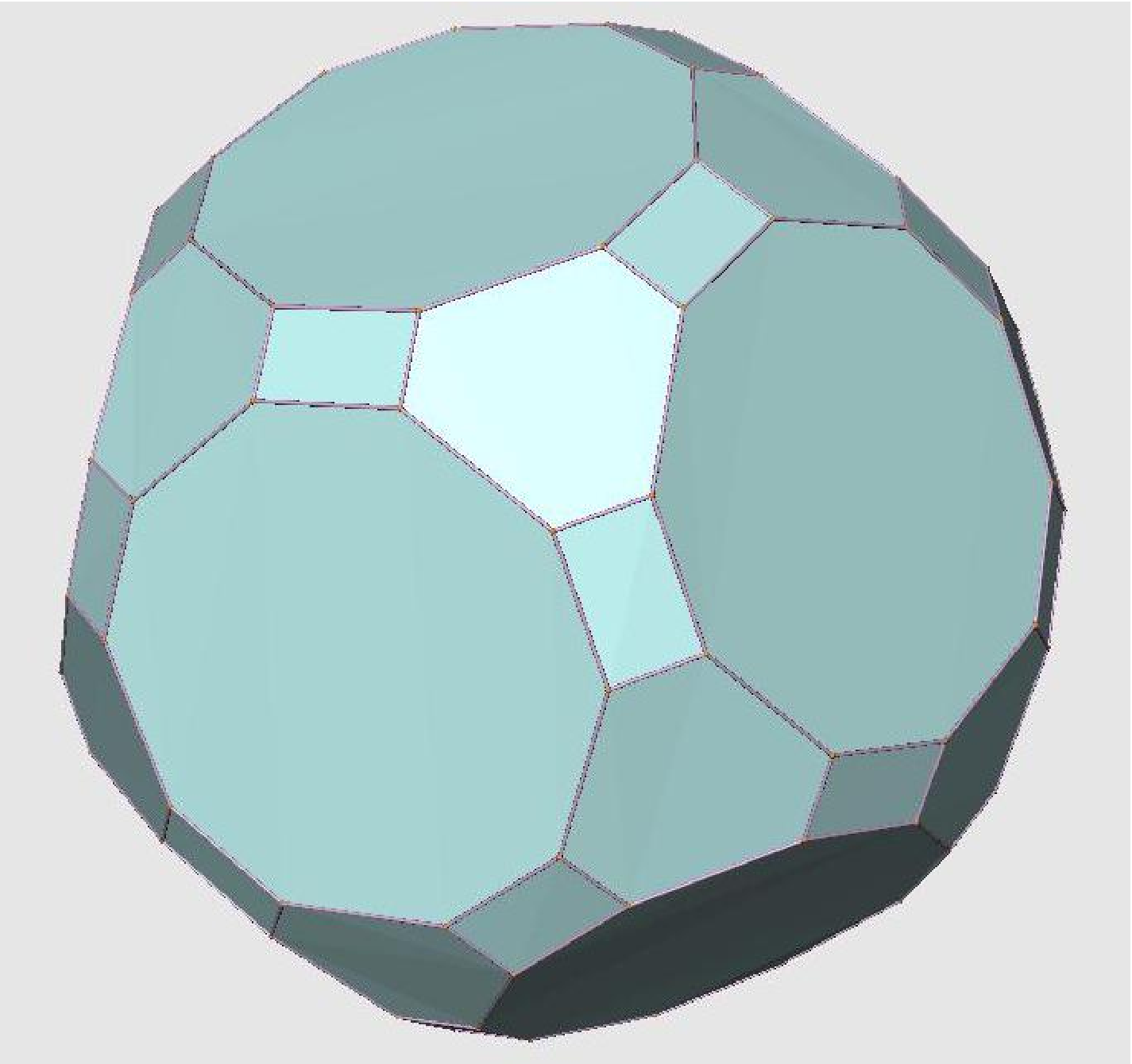}}
  \end{center}
  \caption{The $1$-skeleton of the Archimedean solid $(4,6,10)$ and
    the (non-equilateral) spectral representation of the canonical
    Laplacian for the second highest eigenvalue.}
  \label{fig:specrep111}
\end{figure}

\begin{rmk} \label{rmk:minima} There are analogous versions of Theorem
  \ref{thm:seceigvalarch} in the cases $(4,6,8)$ and $(4,6,6)$. The
  full symmetry group of both solids $(4,6,8)$ and $(4,6,6)$ is the
  full octahedral group, but it is better to view $(4,6,6)$ as a
  polyhedron with the full tetrahedral group (which is a subgroup of
  the full octahedral group) as its symmetry group, by distinguishing
  its hexagonal faces with the help of two colours (say, yellow and
  blue), such that adjacent $6$-gons have different colours. In this
  case the solid $(4,6,6)$ is also called the {\em omnitruncated
    tetrahedron} and has three equivalence classes of edges (separating
  $4$-gons and yellow $6$-gons, $4$-gons and blue $6$-gons, yellow and
  blue $6$-gons), just as the solid $(4,6,8)$ and $(4,6,10)$. The
  corresponding explicit values for $X_0$ and $\lambda_1(X_0)$ are
  $$ X_0 = \frac{1}{13+6\sqrt{2}}(4+\sqrt{2},3+3\sqrt{2},6+2\sqrt{2}) \quad 
  \text{and}\ \lambda_1(X_0) = \frac{11+6\sqrt{2}}{13+6\sqrt{2}} $$ 
  in the case $(4,6,8)$ and
  $$ X_0 = (\frac{3}{10},\frac{3}{10},\frac{2}{5}) \quad 
  \text{and}\ \lambda_1(X_0) = \frac{4}{5}$$
  in the case $(4,6,6)$.
\end{rmk}

Finally, we describe the behaviour of spectral representations
$\Phi_X: V \to S^2$, as $X \in \Delta$ moves towards the boundary
$\partial \Delta$.

\begin{figure}[h]
  \begin{center}
      \psfrag{c1}{$C_1$}
      \psfrag{c2}{$C_2$}
      \psfrag{c3}{$C_3$}
      \psfrag{(4,6,10)}{\small $(4,6,10)$}
      \psfrag{(3,3,3,3,3)}{\small $(3,3,3,3,3)$}
      \psfrag{(3,5,3,5)}{\small $(3,5,3,5)$}
      \psfrag{(3,4,5,4)}{\small $(3,4,5,4)$}
      \psfrag{(3,10,10)}{\small $(3,10,10)$}
      \psfrag{(5,6,6)}{\small $(5,6,6)$}
      \psfrag{(5,5,5)}{\small $(5,5,5)$}
      \psfrag{Simplex}{Simplex $\Delta_\Gamma$}
      \includegraphics[height=4cm]{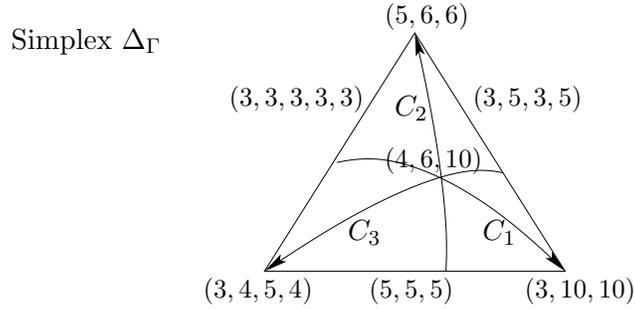}
  \end{center}
  \caption{Convergence behaviour of $\Phi_X$ as $X \to \partial
    \Delta$.}
  \label{fig:archconv1}
\end{figure}

\begin{thm} \label{thm:threecurves} 
  Let $G,\Delta,X_0$ be as in Theorem \ref{thm:seceigvalarch}. Then
  there are three explicitly given curves $C_1,C_2,C_3 \subset \Delta$,
  which meet in $X_0$ and have the following property: For every $X
  \in C_i$, the lengths of two of the three equivalence classes of
  Euclidean edges in the spectral representation of $P_X$ for the
  eigenvalue $\lambda_1(X)$ coincide.

  As $X_n$ converges to the corresponding vertex of the simplex
  $\Delta$ along the curve $C_i$, the spectral representations
  converge to equilateral realisations of the Archimedean solids $(3,10,10)$,
  $(5,6,6)$ and $(3,4,5,4)$, respectively.
  
  For any sequence $X_n \in {\rm int}(\Delta)$ converging to an
  interior point of the boundary edge of the simplex $\Delta$, the
  spectral representations converge to the equilateral realisations of
  one of the solids $(3,3,3,3,3)$, $(5,5,5)$ and $(3,5,3,5)$.

  These convergence properties are illustrated in Figure \ref{fig:archconv1}.
\end{thm}

\begin{figure}[h]
  \begin{center}
  \makebox[4cm]{\includegraphics[width=4cm]{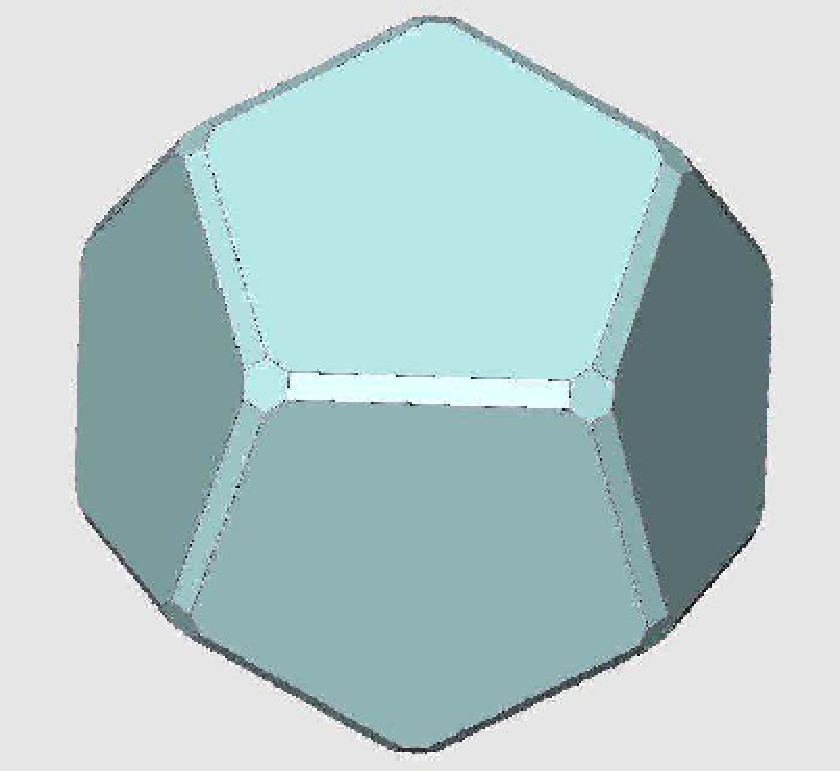}}
  \makebox[4cm]{\includegraphics[width=4cm]{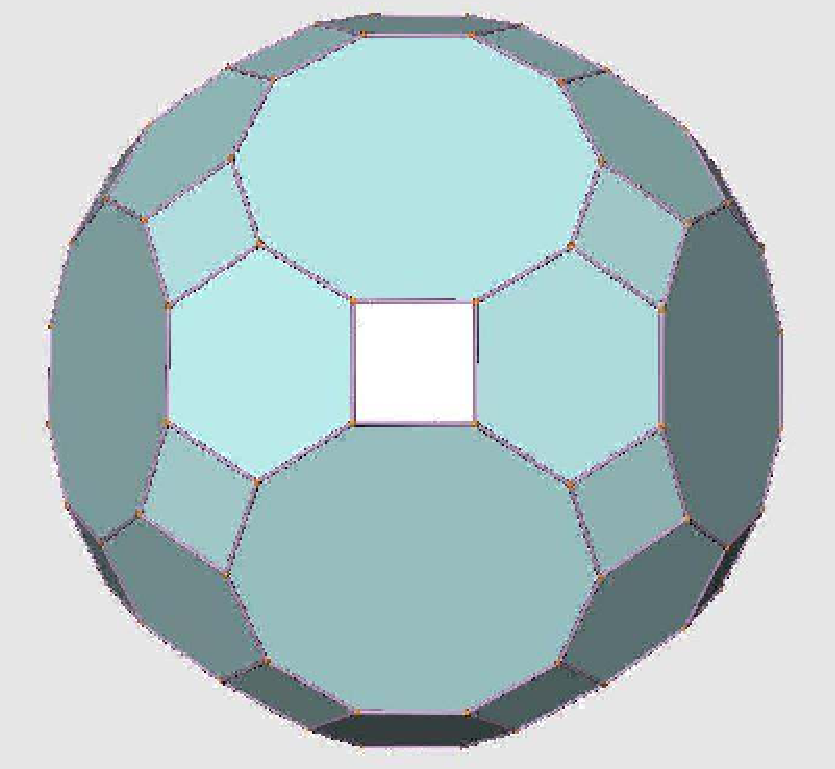}}
  \makebox[4cm]{\includegraphics[width=4cm]{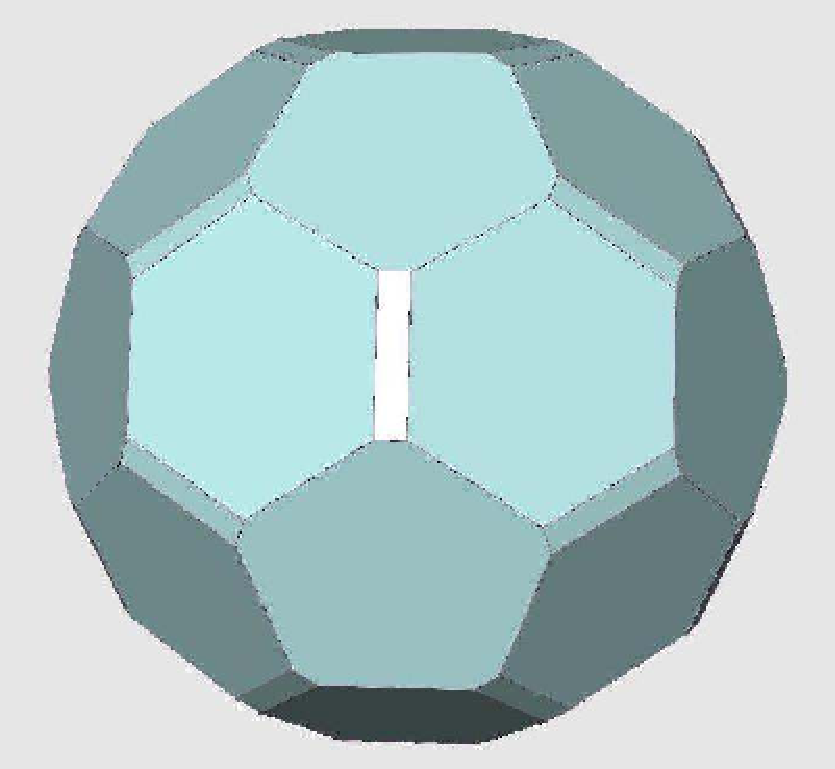}}
  \end{center}
  \caption{Spectral representations of $P_X$ (for points $X$ along $C_2$) 
  for the second highest eigenvalue.}
  \label{fig:transspecrep}
\end{figure}

Figure \ref{fig:transspecrep} shows the spectral representations of
$P_X$ for three points $X$ along the curve $C_2$, illustrating the
transition from the dodecahedron $(5,5,5)$ to the buckeyball
$(5,6,6)$.

\begin{rmk}
  The analogous versions of Theorem \ref{thm:threecurves} for the
  Archimedean solids $(4,6,8)$ and $(4,6,6)$ are illustrated in Figure
  \ref{fig:archconv2} below. The common symmetry group of all solids
  in the diagram containing $(4,6,8)$ is the full octahedral group. In
  the diagram containing $(4,6,6)$, we need to colour the hexagons in
  the solid $(4,6,6)$ with two different colours (as described in
  Remark \ref{rmk:minima}) and, similarly, we have to colour the
  triangles of the solid $(3,4,3,4)$ with two colours such that
  triangles meeting in a vertex have different colours (and refer to
  $(3,4,3,4)$ then as the {\em cantellated tetrahedron}), so that
  the common symmetry group of all solids in this diagram is the full
  tetrahedral group.
  
\begin{figure}[h]
  \begin{center}
      \psfrag{c1}{$C_1$}
      \psfrag{c2}{$C_2$}
      \psfrag{c3}{$C_3$}
      \psfrag{(4,6,8)}{\small $(4,6,8)$}
      \psfrag{(3,3,3,3)}{\small $(3,3,3,3)$}
      \psfrag{(3,4,3,4)}{\small $(3,4,3,4)$}
      \psfrag{(3,4,4,4)}{\small $(3,4,4,4)$}
      \psfrag{(3,8,8)}{\small $(3,8,8)$}
      \psfrag{(4,6,6)}{\small $(4,6,6)$}
      \psfrag{(4,4,4)}{\small $(4,4,4)$}
      \psfrag{(4,6,6)}{\small $(4,6,6)$}
      \psfrag{(3,3,3)}{\small $(3,3,3)$}
      \psfrag{(3,4,3,4)}{\small $(3,4,3,4)$}
      \psfrag{(3,6,6)}{\small $(3,6,6)$}
      \psfrag{(4,6,6)}{\small $(4,6,6)$}
      \includegraphics[height=4cm]{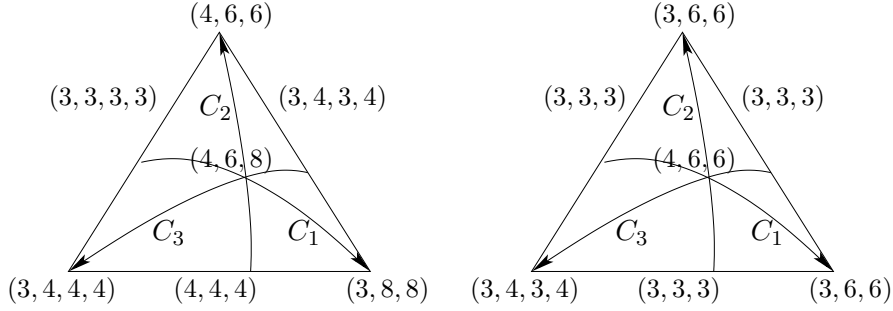}
  \end{center}
  \caption{Boundary convergence behaviour in the cases $(4,6,8)$ and
    $(4,6,6)$.}
  \label{fig:archconv2}
\end{figure}

\end{rmk}

\subsection{Structure of the article}

Section \ref{sec:genresults} provides the proofs of Propositions
\ref{prop:lambda1boundary} and \ref{prop:main2} and of our Main
Theorem \ref{thm:main}. In Sections \ref{sec:coxresults1} and
\ref{sec:coxresults2}, we prove Corollary \ref{cor:coxgroup} and
Propositions \ref{prop:kassabov} and \ref{prop:A3B3H3}. Finally,
Section \ref{sec:arch46x} presents the proofs of Theorems
\ref{thm:seceigvalarch} and \ref{thm:threecurves}.

\medskip

{\bf Acknowlegdements:} We are grateful to Stefan Dantchev, Christoph
Helmberg, Dirk Sch\"utz and Ivan Veseli{\'c} for many helpful
comments. Moreover, we like to express our gratitude to the referee
for providing us with the proof of Proposition \ref{prop:kassabov},
and for many other useful comments which improved and simplifed parts
of our article.

\section{Proofs of the general results} \label{sec:genresults}

Let us start with a convexity proof of $\lambda_1: \Pi_G \to [-1,1]$,
which is used to show that $\lambda_1$ assumes a global minimum in
$\Delta_\Gamma$. 

\smallskip

{\bf Proof of Proposition \ref{prop:lambda1boundary}:} 
First note that, for $P, P' \in \Pi_G$ and $\alpha \in [0,1]$,
$$ \frac{\left\langle (\alpha P + (1-\alpha)P')f,f \right\rangle}
{\Vert f \Vert^2} =
\alpha \frac{\langle Pf, f \rangle}{\Vert f \Vert^2} + (1-\alpha)
\frac{\langle P' f,f \rangle}{\Vert f \Vert^2}, $$ 
which implies the convexity of $\lambda_1: \Pi_G \to [-1,1]$ by taking
supremums on both sides and using the characterisation
\eqref{eq:rayleigh}.

Note also that $\Pi_G$ is compact, and the continuous function
$\lambda_1: \Pi_G \to [-1,1]$ must have a global minimum at some point
$P \in \Pi_G$. If $P \not\in \Delta_\Gamma$, then consider the
$\Gamma$-invariant matrix $P' = \frac{1}{|\Gamma|} \sum_{\gamma \in
  \Gamma} \gamma P \in \Pi_G$, and we conclude $\lambda_1(P') \le
\lambda_1(P)$ by the convexity of $\lambda_1$. $P'$ may have a
non-vanishing (constant) main diagonal. Nevertheless, we can write $P'
= \beta {\rm Id} + (1-\beta) P''$ with appropriate $P'' \in
\Delta_\Gamma$ and $\beta \ge 0$. This implies that
$$ \lambda_1(P') - \lambda_1(P'') = \beta (1- \lambda_1(P'')) \ge 0, $$
which shows that $\lambda_1$ assumes also a global minimum at $P'' \in
\Delta_\Gamma$.

Now, assume that $G = \Cay(\Gamma,S)$ and that $S$ is a minimal set of
generators. The minimality of $S$ implies that, for every $X
\in \partial \Delta_\Gamma$, the graph $G_{P_X}$ consists of more than
one connected component and, therefore, $\lambda_0(P_X) =
\lambda_1(P_X) = 1$. Then \eqref{eq:convbound} follows from the
continuity of the function $\lambda_1$. For every interior point $X
\in \Delta_\Gamma$ we have $\lambda_1(P_X) < 1$, since $G_{P_X}$ is
connected.  \EPf

Our next goal is the proof that $\Gamma$-invariant spectral representations 
map all vertices onto a sphere and that equivalent edges are mapped to
Euclidean segments of the same length.

\begin{lem} \label{lem:basic} 
  Let $G=(V,E)$ be a simple finite graph and $\Gamma \subset \Aut(G)$
  be vertex transitive with $N$ equivalence classes of edges, and
  $\Delta_\Gamma$ be as in \eqref{eq:deltagen}. Let $X \in \Delta_\Gamma$
  and $\lambda \in [-,1,1]$ be an eigenvalue of multiplicity $k \ge 2$
  of the operator $P_X$. Let $\Phi=\Phi_{X,\lambda}: V \to \R^k$ be
  the associated spectral representation.  Then there exist constants
  $c >0$, $c_1,\dots,c_N \ge 0$ such that
  \begin{itemize}
  \item[(a)] for all vertices $i \in V$: $\Vert \Phi(i) \Vert = c$,
  \item[(b)] for all edges $\{i,j\} \in E$ in the $l$-th equivalence class
  $$ \Vert \Phi(i) - \Phi(j) \Vert = c_l. $$
  \end{itemize}
\end{lem}

{\bf Proof:} Let $\phi_1,\dots,\phi_k$ be the orthonormal basis of
eigenfunctions defining $\Phi = (\phi_1,\dots,\phi_k)^\top$ (we
consider $\Phi(i)$ as a column vector). Let $\gamma \in \Gamma$ be
fixed and $\psi_r = \phi_r \circ \gamma: V \to \R$. One easily checks
that $\psi_1,\dots,\psi_k$ are also an orthonormal basis satisfying
$P_X \psi_r = \lambda \psi_r$. Consequently, there exists a matrix
$C = (c_{rs}) \in O(k)$ such that $\psi_r = \sum_{s=1}^k c_{rs} \phi_s$.
This implies $\Phi(\gamma i) = C\Phi(i)$ and
\begin{equation} \label{eq:relgammaij} 
  \langle \Phi(\gamma i),\Phi(\gamma j) \rangle = 
  \langle C\Phi(i),C\Phi(j) \rangle= 
  \langle \Phi(i),\Phi(j) \rangle.
\end{equation}
\eqref{eq:relgammaij} implies (a) by choosing $i=j$ and using the
vertex transitivity of $\Gamma$. (b) follows from (a),
\eqref{eq:relgammaij} and
$$ \Vert \Phi(i) - \Phi(j) \Vert^2 = \Vert \Phi(i) \Vert^2 - 
2 \langle \Phi(i),\Phi(j) \rangle + \Vert \Phi(j) \Vert^2. $$ \EPf

Before entering into the proof of our main result, let us remark that
the above identity \eqref{eq:relgammaij} can be rewritten as
\begin{equation} \label{eq:gammainv}
  \sum_{r=1}^k \phi_r(i) \phi_r(j) = 
  \sum_{r=1}^k \phi_r(\gamma i) \phi_r(\gamma j) \quad 
  \text{for all $i,j \in V$ and $\gamma \in \Gamma$}. 
\end{equation}
Moreover, observe that the following identity is an immediate
consequence of the vertex transitivity of $\Gamma$, the left coset
decomposition $\Gamma = \gamma_1 \Gamma_1 \cup \gamma_2 \Gamma_1 \cup
\cdots \cup \gamma_n \Gamma_1$, where $\Gamma_1 \subset \Gamma$ is the
stabilizer of $1 \in V$, and the relation $|\Gamma| = n |\Gamma_1|$:
\begin{equation} \label{eq:fundeq}
  \sum_{i=1}^n f(i) = \frac{n}{|\Gamma|} \sum_{\gamma \in \Gamma} f(\gamma 1).
\end{equation}

{\bf Proof of Theorem \ref{thm:main}:} 
For simplicity, we discuss the key arguments for the special choice of
the first and second equivalence class of edges. The proof for two
arbitrary equivalence classes is completely analogous.

Let $\{1,i_1\},\dots,\{1,i_p\}$ be all edges adjacent to $1 \in V$ in
the first equivalence class of edges (note that $p=m_1$). Let
$\{1,j_1\},\dots,\{1,i_q\}$ be all edges adjacent to $1 \in V$ in the
second equivalence class of edges (note that $q=m_2$).  We conclude
from \eqref{eq:gammainv} that
\begin{equation} \label{eq:i1ip}
\sum_{r=1}^k \phi_r(1) \phi_r(i_1) = 
   \sum_{r=1}^k \phi_r(1) \phi_r(i_2) = \cdots = 
   \sum_{r=1}^k \phi_r(1) \phi_r(i_p), 
\end{equation}
and the same identity holds for the edges in the second equivalence class.

Let $X_0 \in U$ an arbitary point (not necessarily a critical point of
$\lambda_1$), $\xi = (\frac{1}{p},-\frac{1}{q},0,\dots,0) \in \R^N$
and $X_t := X_0 + t \xi \in U$ for $t \in (-\epsilon,\epsilon)$ and
$\epsilon > 0$ suitably small. For simplicity of notation, we
introduce $P(t) := P_{X_t}$, $\lambda(t)=\lambda(X_t)$. Let
$\phi_1,\dots,\phi_k$ be an orthonormal basis of the eigenspace
$\E_{\lambda(0)}(P(0))$. Let ${\rm Pr}_t$ denote the orthogonal
projection of $l^2(G)$ onto the eigenspace $\E_{\lambda(t)}(P(t))={\rm
  ker}(P(t)-\lambda(t))$. Since $P(t)$ and $\lambda(t)$ depend
smoothly on $t$, ${\rm Pr}_t$ is also smooth in $t$. By making
$\epsilon > 0$ smaller, if needed, we can assume that ${\rm
  Pr}_t\phi_1,\dots,{\rm Pr}_t\phi_k$ is a basis of
$\E_{\lambda(t)}(P(t))$, for all $t \in
(-\epsilon,\epsilon)$. Applying Gram-Schmidt to these vectors, we
obtain an orthonormal basis $\phi_{1,t}, \dots, \phi_{k,t}$ of the
eigenspace $\E_{\lambda(t)}(P(t))$, depending smoothly on $t$ and
satisfying $\phi_r = \phi_{r,0}$. Note that $P'(0) = (c_{ij})$ with
$$ c_{\gamma 1,\gamma i} = \begin{cases} 
\frac{1}{p}, & \text{if $i \in \{i_1,\dots,i_p\}$}, \\ 
- \frac{1}{q}, & \text{if $i \in \{j_1,\dots,j_q\}$}, \\ 
0, & \text{otherwise}, \end{cases} $$
for all $\gamma \in \Gamma$ and $i \in V$.  

Let $r \in \{1,\dots,k\}$. By the orthonormality of the functions
$\phi_r$, we have $\langle \phi_r, \frac{\partial}{\partial
  t}\big\vert_{t=0} \phi_{r,t} \rangle = 0$. Using this and the
symmetry of $P(t)$, we obtain, by differentiating $\lambda(t) =
\langle P(t) \phi_{r,t},\phi_{r,t} \rangle$ at $t=0$:
\begin{eqnarray*}
  \lambda'(0) = \langle P'(0) \phi_r, \phi_r \rangle &=& 
  \sum_{i=1}^n \phi_r(i) \left( \sum_{j=1}^n c_{ij} \phi_r(j) \right) \\
  &=& \frac{n}{|\Gamma|} \sum_{\gamma \in \Gamma} \phi_r(\gamma 1) 
  \left( \sum_{j=1}^n c_{\gamma 1,\gamma j} \phi_r(\gamma j) \right) 
  \qquad \text{using \eqref{eq:fundeq}} \\
  &=& \frac{n}{|\Gamma|} \sum_{\gamma \in \Gamma} \phi_r(\gamma 1) 
  \left( \frac{1}{p} \sum_{s=1}^p \phi_r(\gamma i_s) - 
    \frac{1}{q} \sum_{s=1}^q \phi_r(\gamma j_s) \right).
\end{eqnarray*}
On the other hand, we have
\begin{eqnarray*}
\langle \Phi(1),\Phi(i_1) \rangle &=& 
\frac{1}{p} \sum_{s=1}^p \sum_{r=1}^k \phi_r(1) \phi_r(i_s) 
\qquad \text{using \eqref{eq:i1ip}} \\
&=& \sum_{r=1}^k \frac{1}{|\Gamma|} \sum_{\gamma \in \Gamma} 
\phi_r(\gamma 1) \frac{1}{p} \sum_{s=1}^p \phi_r(\gamma i_s)  
\qquad \text{using \eqref{eq:gammainv}}.
\end{eqnarray*}
Combining both results, we obtain
$$ \langle \Phi(1),\Phi(i_1) \rangle - \langle \Phi(1),\Phi(j_1) \rangle 
= \frac{k}{n} \lambda'(0). $$ 
This implies that we have $\Vert \Phi(1) - \Phi(i_1) \Vert = \Vert
\Phi(1) - \Phi(j_1) \Vert$ if and only if $\lambda'(0)=0$. 

Since the above arguments hold for any choice of equivalence classes
of edges, we have $c_1=\dots=c_N$ in Lemma \ref{lem:basic} above
(i.e., an equilateral spectral representation) if and only if the
derivative of $\lambda$ at $X_0$ vanishes in all directions of the
simplex, i.e., if $X_0 \in U$ is a critical point of $\lambda$.  \EPf

The proof of Proposition \ref{prop:main2} is based on the following lemma:

\begin{lem} \label{lem:relinnerprod} Let $\Gamma$ be a finite group,
  $\rho: \Gamma \to O(k)$ be an irreducible representation and, as
  before, $\langle \cdot, \cdot \rangle$ be the standard inner product
  in $\R^k$. For any non-zero vector $p \in \R^k$ there is a constant
  $C_p > 0$ such that
  $$ \sum_{\gamma \in \Gamma} \langle \rho(\gamma) p,v \rangle \, \langle
  \rho(\gamma) p,w \rangle = C_p \langle v,w \rangle $$
  for all $v,w \in \R^k$.
\end{lem}

{\bf Proof:} The expression
  $$ \langle v,w \rangle_p := \sum_{\gamma \in \Gamma} \langle \rho(\gamma) p,v 
  \rangle \, \langle \rho(\gamma) p,w \rangle $$ is obviously a symmetric
  bilinear form. The form is positive definite, since $\langle v, v
  \rangle_p = 0$ implies that $v$ is perpendicular (w.r.t. the
  standard inner product) to $\spann \{\rho(\gamma) p \mid \gamma \in
  \Gamma \}$. Irreducibility of $\rho$ implies that $\spann
  \{\rho(\gamma) p \mid \gamma \in \Gamma \}= \R^k$, so $v =
  0$. Therefore, there exists a positive definite symmetric matrix $A$
  such that
  $$ \langle v,w \rangle_p = \langle Av,w \rangle. $$ 
  Let $\gamma_0 \in \Gamma$. Then
  $$ \langle \rho(\gamma_0)v,\rho(\gamma_0)w \rangle_p = \sum_{\gamma \in \Gamma} 
  \langle \rho(\gamma_0^{-1}\gamma) p, v \rangle \langle 
  \rho(\gamma_0^{-1}\gamma) p, w \rangle = \langle v,w \rangle_p, $$ 
  i.e., $\langle \cdot, \cdot \rangle_p$ is $\rho(\Gamma)$-invariant,
  and we have $A \rho(\gamma) = \rho(\gamma) A$ for all $\gamma \in
  \Gamma$. Since $\rho$ is irreducible, we conclude from Schur's lemma
  that $A$ is of the form $C_p \cdot \Id$ with a constant $C_p >
  0$. This finishes the proof of the lemma.  \EPf

{\bf Proof of Proposition \ref{prop:main2}:} Note that the vertices of
$G = \Cay(\Gamma,S)$ are the group elements, and that
$$ P_X f(\gamma) = \sum_{j=1}^\nu x_j f(\gamma s_j) + 
\sum_{j=1}^\mu x_{\nu+j} (f(\gamma \sigma_j)+f(\gamma \sigma_j^{-1})). $$
This implies that
\begin{eqnarray*} 
P_X \phi_r(\gamma) &=& \sum_{j=1}^\nu x_j \pi_r(\rho(\gamma s_j)p) 
+ \sum_{j=1}^\mu x_{\nu+j} ( \pi_r(\rho(\gamma \sigma_j)p) + 
\pi_r(\rho(\gamma \sigma_j^{-1})p) ) \\
&=& \pi_r\left( \rho(\gamma) \left( \sum_{j=1}^\nu x_j \rho(s_j)p + 
\sum_{j=1}^\mu x_{\nu+j} ( \rho(\sigma_j)p + \rho(\sigma_j^{-1})p ) 
\right) \right) \\
&=& \lambda \pi_r(\rho(\gamma)p) = \lambda \phi_r(\gamma),
\end{eqnarray*}
by using \eqref{eq:lamP}. This shows that $\phi_r$ is an eigenfunction
of $P_X$ for the eigenvalue $\lambda$. The orthogonality of the
functions $\phi_r$ is a straightforward application of Lemma
\ref{lem:relinnerprod}:
$$
\langle \phi_r,\phi_s \rangle = \sum_{\gamma \in \Gamma}
\pi_r(\rho(\gamma)p) \pi_s(\rho(\gamma)p) = \sum_{\gamma \in \Gamma}
\langle \rho(\gamma)p,e_r \rangle \, \langle \rho(\gamma)p,e_s
\rangle = C_P \langle e_r, e_s \rangle,
$$
where $e_1,\dots,e_k$ denotes the standard basis in $\R^k$. Now let
$\Gamma = \{ \gamma_1,\dots,\gamma_n \}$. Let $A$ be the $(k \times
n)$ matrix whose columns are the vectors $\rho(\gamma_j)p \in
S^{k-1}$. Then the rows of $A$ represent the functions $\phi_r$, and
we have
$$ \sum_{r=1}^k \Vert \phi_r \Vert^2 = 
\sum_{j=1}^n \Vert \rho(\gamma_j)p \Vert^2 = n = | \Gamma |. $$
This shows that $\Vert \phi_r \Vert^2 = \frac{|\Gamma|}{k}$. \EPf

\begin{rmk}
  Assume that $\rho$ in Proposition \ref{prop:main2} is irreducible
  but not orthogonal, i.e., $\rho: \Gamma \to GL(k,\R)$. The above
  proof still shows that the functions $\phi_r$ are
  eigenfunctions. Let $A$ be the $(k \times n)$ matrix as in the
  proof. Then the irreducibility of $\rho$ implies that the columns
  span all of $\R^k$, i.e., the rank of $A$ is $k$. But this means that
  the functions $\phi_r$ (the $k$ rows of $A$) must be linearly independent.
\end{rmk}

\section{Proof of Corollary \ref{cor:coxgroup} and Proposition
  \ref{prop:kassabov}} \label{sec:coxresults1}

Our first aim is to establish the geometric procedure for obtaining
eigenfunctions of $P_X$ on the Cayleygraph of a Coxeter group, as well
as explicit derivations of the maps $\Psi_\lambda$ and $\Psi_\Delta$.

\smallskip

{\bf Proof of Corollary \ref{cor:coxgroup}:} We start with a finite,
irreducible Coxeter group. This implies that the geometric realisation
$\Gamma \hookrightarrow O(k)$ is an irreducible, faithful
representation. Note that we have $\langle n_i,n_j \rangle = -\cos
\frac{\pi}{m_{ij}} $, where $m_{ij}$ is the order of the element $s_i
s_j$. Since $\Gamma$ is a finite Coxeter group, $M = ( \langle n_i,n_j
\rangle)$ is a positive definite, symmetric matrix. Writing $M = \Id
-C$ with a symmetric matrix $C < {\rm Id}$ (as quadratic forms) whose
entries are all non-negative, we obtain $M^{-1} = \sum_{s=0}^\infty
C^s$. Irreducibility implies that for every position $1 \le i,j \le
k$, there is an $s \ge 0$ such that $(C^s)_{ij} > 0$. This implies
that all entries of $M^{-1}$ are strictly positive. Recall that
$$ V := \det(n_1,\dots,n_k) = \left( \det M \right)^{1/2} > 0. $$
We define, as in \eqref{eq:Pj},
$$ p_j = (-1)^{j-1} n_1 \times \cdots \times \widehat{n_j} \times \cdots
\times \cdots n_k. $$ 

The vectors $p_j$ may all have different Euclidean lengths. We have,
by construction $\langle n_i,p_j \rangle = V \delta_{ij}$. Let
$$ 
\Delta = \Delta_\Gamma = \{ (x_1,\dots,x_k) \mid x_j \ge 0, \sum_j x_j = 1 \} 
$$
be the simplex associated to the Cayley graph $\Cay(\Gamma,S)$. 

Our aim is to construct the maps $\Psi_\Delta: \F_0 \to {\rm int}(\Delta)$ and
$\Psi_\lambda: \F_0 \to (0,1)$: Any point $p \in \F_0$ can be
expressed uniquely as
$$ p = \alpha_1 p_1 + \cdots + \alpha_k p_k, $$
with $\alpha_1,\dots,\alpha_k > 0$. We will show that there is a
unique choice of $X=(x_1,\dots,x_k) \in {\rm int}(\Delta)$ and
$\lambda \in (0,1)$ such that
\begin{equation} \label{eq:eigeqP}
\lambda p = \sum_j x_j \sigma_j(p). 
\end{equation}
We then define $\Psi_\Delta(p) = X$ and $\Psi_\lambda(p) =
\lambda$. The construction will show that $X$ and $\lambda$ depend
smoothly on the coordinates $\alpha_j$. Applying Proposition
\ref{prop:main2} yields the results stated in the Corollary. It then
only remains to prove that $\Psi_\Delta$ is bijective and that the
composition $\Psi_\lambda \circ \Psi_\Delta^{-1}$ is analytic.

Since $\sigma_j(p) = p - 2 \langle p,n_j \rangle n_j = p - 2 \alpha_j V n_j$,
we immediately see that $\Vert p - \sigma_j(p) \Vert = 2 \alpha_j V$. Moreover,
\eqref{eq:eigeqP} translates into
\begin{equation} \label{eq:eigeqP2}
\lambda p = (x_1 + \cdots + x_k)p - 2 V \sum_j \alpha_j x_j n_j. 
\end{equation}
This means that we need to find a unique $(x_1,\dots,x_n) \in \Delta$
and $\mu \in \R$ such that
\begin{equation} \label{eq:nPeq}
\sum_j \alpha_j x_j n_j = \mu \sum_j \alpha_j p_j,
\end{equation}
and then set $\lambda = x_1+\cdots+x_k-2V\mu$. Taking inner products
with the simple roots $n_1, \dots, n_k$, and bringing everything in a
matrix equation, we end up with the equivalent equation
\begin{equation} \label{eq:crucialax} 
M \begin{pmatrix} \alpha_1 x_1 \\ \vdots \\ \alpha_k x_k \end{pmatrix}
= \mu V \begin{pmatrix} \alpha_1 \\ \vdots \\ \alpha_k \end{pmatrix}. 
\end{equation}
Obviously, this equation is homogeneous, i.e., if
$(x_1,\dots,x_k,\mu)$ is a solution then so is $(c x_1,\dots,c x_k,c
\mu)$ for any constant $c$. We first seek for the unique solution
$(x_1',\dots,x_k')$ of \eqref{eq:crucialax} for the choice
$\mu=1$. $(x_1',\dots,x_k')$ will not be a point in $\Delta$, and we
obtain the correct solution by way of rescaling. Using the fact that
$M^{-1} = \Id + D$, where all diagonal entries of $D$ are non-negative
and all off-diagonal are strictly positive, we end up with the
inequality
\begin{equation} \label{eq:calcX'}
\begin{pmatrix} x_1' \\ \vdots \\ x_k' \end{pmatrix} =
V \begin{pmatrix} \alpha_1 & & \\ & \ddots & \\ & &
  \alpha_k \end{pmatrix}^{-1} M^{-1} \begin{pmatrix} \alpha_1 & & \\ &
  \ddots & \\ & & \alpha_k \end{pmatrix} \begin{pmatrix} 1 \\
  \vdots \\ 1 \end{pmatrix} > V \begin{pmatrix} 1 \\ \vdots \\
  1 \end{pmatrix}.
\end{equation}
This shows that any choice $\alpha_1,\dots,\alpha_k > 0$ leads to
a strictly positive vector $(x_1',\dots,x_k')$, and that
\begin{equation} \label{eq:calcl'}
\lambda' := x_1' + \ + \cdots + x_k' - 2 V > 0. 
\end{equation}
For $p = \sum_j \alpha_j p_j \in \F_0$, we first calculate
$x_1',\dots,x_k',\lambda'> 0$ via the equations \eqref{eq:calcX'} and
\eqref{eq:calcl'}, and then apply the rescaling to obtain
\begin{equation} \label{eq:PsisP}
\Psi_\Delta(p) = \frac{1}{\sum_j x_j'}(x_1',\dots,x_k') \in {\rm int}(\Delta), 
\qquad \Psi_\lambda(p) = \frac{\lambda'}{\sum_j x_j'} \in (0,1).
\end{equation}
 
Next we show that $\Psi_\Delta: \F_0 \to {\rm int}(\Delta)$ is
bijective. Choose $X=(x_1,\dots,x_k) \in {\rm int}\Delta$. An equivalent
reformulation of \eqref{eq:crucialax} is
\begin{equation} \label{eq:eigequation}
  \begin{pmatrix} \alpha_1 \\ \vdots \\ \alpha_k \end{pmatrix} = \mu
  V {\mathbf x}^{-1} M^{-1} \begin{pmatrix} \alpha_1 \\ \vdots \\ \alpha_k 
  \end{pmatrix},
\end{equation}
where ${\mathbf x} = \diag(x_1,\dots,x_k)$ denotes the diagonal matrix
with the entries $x_j$. Note that $V\, {\mathbf x}^{-1} M^{-1}$ is a
matrix with all its entries strictly positive. Therefore, we can apply
Perron-Frobenius theory and conclude that there is a unique
Perron-Frobenius eigenvector $(\alpha_1,\dots,\alpha_k)$, scaled in
such a way that $p = \sum_j \alpha_j p_j \in S^{k-1}$. Since $\alpha_j
> 0$, we conclude that $p \in \F_0$.  This shows that every $X \in
{\rm int}(\Delta)$ has a unique preimage under $\Psi_\Delta$.

Moreover, note that $\mu^{-1}$ is the Perron-Frobenius eigenvalue of
the matrix $V {\mathbf x}^{-1} M^{-1}$ and that $\lambda =
\left(\sum_j x_j\right) - 2 V \mu = 1 - 2 V \mu$ in
\eqref{eq:eigeqP2}. This implies that the composition $\Psi_\lambda \circ
\Psi_\Delta^{-1}: {\rm int}(\Delta) \to (0,1)$ is given by
\begin{equation} \label{eq:PsiPsi}
\Psi_\lambda \circ \Psi_\Delta^{-1}(X) = 1 -  2 \Lambda_X, 
\end{equation}
where $\Lambda_X$ is the Perron-Frobenius eigenvalue of the positive
matrix ${\mathbf x}^{-1} M^{-1}$. Since this eigenvalue has always
multiplicity one, it depends analytically on the weights
$x_1,\dots,x_k$, by the analytic version of the Implicit Function
Theorem. This finishes the proof of Corollary \ref{cor:coxgroup}. \EPf

\smallskip

Next we modify arguments in Kassabov \cite[p. 20]{Kas09} to prove
Proposition \ref{prop:kassabov}. 

\smallskip

{\bf Proof of Proposition \ref{prop:kassabov}:} Let us first recall
some of his notation of this source. Let ${\mathcal H} = l^2(G)$ and
$\pi: \Gamma \to U({\mathcal H})$ be the right-regular representation
$(\pi(\sigma) f)(\gamma) = f(\gamma \sigma)$. Let
$$ \Delta_X = {\rm Id} - P_X = \sum_{i=1}^k x_i ({\rm Id} - \pi(s_i)). $$
Let $V_i = \{f \in {\mathcal H} \mid f = \pi(s_i) f \}$. Note that
  ${\rm Id} - \pi(s_i)$ is equal to $2 {\rm Pr}_{V_i^\bot}$, i.e.,
  twice the orthogonal projection to the orthogonal complement of the
  subspace $V_i$. This implies that
\begin{equation} \label{eq:kas1} 
\langle \Delta_X f,f \rangle = 2 \sum_i x_i d_{V_i}(f)^2, 
\end{equation}
where $d_V(f) = \inf_{g \in V} \Vert f-g \Vert = \Vert {\rm
  Pr}_{V^\bot} f \Vert$. Moreover, we have
$$ \bigcap_{i=1}^k V_i = \{ \text{constant functions in $l^2(G)$} \}. $$
For $f \in {\mathbb H}$, let ${\mathbf d}_f$ denote the column vector
with entries the distances $d_{V_i}(f)$. We conclude from
\cite[Thm. 5.1]{Kas09} that, for any function $f$ orthogonal to the
constant functions,
\begin{equation} \label{eq:kas2} \Vert f \Vert^2 \le {\mathbf
    d}_f^\top M^{-1} {\mathbf d}_f = ({\mathbf x}^{1/2} {\mathbf
    d}_f)^\top {\mathbf x}^{-1/2} M^{-1} {\mathbf x}^{-1/2} ({\mathbf
    x}^{1/2} {\mathbf d}_f) \le \Lambda_X \Vert {\mathbf x}^{1/2}
  {\mathbf d}_f \Vert^2,
\end{equation}
where $\Lambda_X$ is the the Perron-Frobenius eigenvalue of ${\mathbf
  x}^{-1/2} M^{-1} {\mathbf x}^{-1/2}$ (note that this agrees with the
Perron-Frobenius eigenvalue of ${\mathbf x}^{-1} M^{-1}$). Combining
\eqref{eq:kas1} and \eqref{eq:kas2}, we conclude that
$$ \langle \Delta_X f, f \rangle \ge 2 \Lambda_X \Vert f \Vert^2, $$
i.e., the second highest eigenvalue $\lambda_1(X)$ of $P_X$ is $\le 1-
2 \Lambda_X$. On the other hand, \eqref{eq:PsiPsi} in the previous
proof shows that $\Psi_\lambda \circ \Psi_\Delta^{-1}(X) = 1-2
\Lambda_X$ is a non-trivial eigenvalue of $P_X$ (of multiplicity $\ge
k$), and therefore we must have $\lambda_1 = \Psi_\lambda \circ
\Psi_\Delta^{-1}$. This finishes the proof of Proposition
\ref{prop:kassabov}. \EPf

\section{Proof of Proposition \ref{prop:A3B3H3}} 
\label{sec:coxresults2}

Our main goal is to prove Corollary \ref{cor:mult3} below. We follow
closely the arguments given by van der Holst \cite{vdH95}. We use the
notation used there, but recall them for the reader's convenience. Let
$G=(V,E)$ be an arbitrary connected graph with vertex set $V =
\{1,\dots,n\}$. For a given subset $V_0 \subset V$ of vertices, we
define $\langle V_0 \rangle \subset G$ to be the subgraph induced by
$V_0$. For a function $f \in l^2(G)$, let $\supp(f) := \{ i \in V \mid
f(i) \neq 0 \}$ and $\supp_{\pm}(f) = \{ i \in \supp(f) \mid \pm f(i)
> 0 \}$. We say that a non-zero function $f$ in a subspace $\E \subset
l^2(G)$ has {\em minimal support}, if for every non-zero function $g
\in \E$ with $\supp(g) \subset \supp(f)$ we have $\supp(g) =
\supp(f)$.

Let ${\mathcal M}(G)$ be the set of all symmetric (not necessarily
stochastic) matrices $M = (m_{ij})$ with all non-diagonal entries
$m_{ij} > 0$ if $i \sim j$ and $m_{ij} = 0$ if $i \not\sim j$.  Note
that we do not impose any sign conditions on the diagonal entries
$m_{ii}$. It is a direct consequence of the connectedness of $G$ and
Perron-Frobenius that the highest eigenvalue $\lambda_0(M)$ is
simple. Colin de Verdi{\'e}re calls the matrices in ${\mathcal M}(G)$
{\em Schr\"odinger operators on the graph $G$} (see \cite{CdV98}), and
they play an important role for his graph invariant $\mu(G)$ (see
\cite{CdV91}). The following result can be considered as a graph
theoretical analogue of the Courant nodal domain for Riemannian
manifolds (see \cite{Cha84}):

\begin{prop}[\cite{vdH95}] \label{prop:courant} Let $G = (V,E)$ be a
  finite connected graph and $M \in {\mathcal M}(G)$. Let
  $\E=\E_{\lambda_1(M)}$ be the eigenspace of the second highest
  eigenvalue of $M$. Let $f \in \E$ be a function of minimal
  support. Then $\langle \supp_+(f) \rangle$ and $\langle \supp_-(f)
  \rangle$ are both connected graphs.
\end{prop}

This fact allows us to prove the following special result:

\begin{cor} \label{cor:mult3} 
  Let $\Gamma \in \{ A_3, B_3, H_3 \}$ and $G$ be the associated
  Cayley graph with respect to the canonical set $S = \{s_1,s_2,s_3\}$
  of generators. Let $X$ be an interior point of $\Delta_\Gamma$. Then
  we have $\lambda_1 = \lambda_1(X) \in [0,1)$, and the corresponding
  eigenspace has dimension $\le 3$.
\end{cor}

{\bf Proof (following mainly \cite{vdH95}):} Let $\E$ be the
eigenspace of $\lambda_1$. Since $X \in {\rm int}(\Delta)$, we have
$P_X \in {\mathcal M}(G)$. Note that the spectrum of $P_X$ is
symmetric (since $G$ is bipartite), and therefore we must have
$\lambda_1 \in [0,1)$, since both eigenvalues $-1,1$ are simple,
because $G_{P_X}=G$ is connected.

Recall that $V = \Gamma$ and that $G = (V,E)$ is the one-skeleton of
one of the solids $(4,6,6)$, $(4,6,8)$ or $(4,6,10)$. In particular,
$G$ is a $3$-connected finite planar graph of constant vertex degree
three. We think of the elements of $G$ as being enumerated and
identify group elements with their corresponding integers. Thus, it
makes sense to write $p_{\gamma,\gamma'}$ for the matrix entries of
$P_X$.

Let $\gamma_0 \in V$ be arbitrary and $\Theta_{\gamma_0}: \E \to \R^3$
be the map
\begin{equation} \label{eq:theta} 
  \Theta_{\gamma_0}(f) = (f(\gamma_0 s_1),f(\gamma_0 s_2),f(\gamma_0 s_3)).
\end{equation}
We prove that this map is {\em injective}, which shows that ${\rm dim}\, \E
\le 3$. Assume that there is a non-zero $f \in \E$ with $\Theta(f) = 0$,
i.e., $\supp(f) \cap \gamma_0 S = \emptyset$. Choose a function $g \in
\E$ with minimal support $\supp(g) \subset \supp(f)$.

We first show that $g(\gamma_0) = 0$. Assume that $g(\gamma_0) \neq
0$. Without loss of generality, we can assume that $\gamma_0 \in
\supp_+(g)$ (otherwise replace $g$ by $-g$).  Since
$$ 
\lambda_1 g(\gamma_0) = \sum_{j=1}^3 p_{\gamma_0,\gamma_0 s_j} g(\gamma_0 s_j) = 0, 
$$  
we must have $\lambda_1=0$. Since $\supp_+(g)$ is connected by
Proposition \ref{prop:courant} and $g$ vanishes on all neighbours of
$\gamma_0$, we conclude that $\supp_+(g) = \{ \gamma_0 \}$. Let
$S_n(\gamma) \subset V$ denote the sphere of combinatorial radius $n$
around $\gamma$. Since for our graphs, all vertices in $S_1(\gamma_0)$
have two neighbours in $S_2(\gamma_0)$ and $g$ is an
eigenfunction to the eigenvalue zero, we must have $f(\gamma') \le 0$
for all $\gamma' \in S_2(\gamma_0)$, and there exists a $\gamma_1 \in
S_2(\gamma_0)$ with $f(\gamma_1) < 0$. Now, $\gamma_1$ cannot be a
neighbour of all three vertices in $S_1(\gamma_0)$, and therefore must
have a neighbour $\gamma_2$ with distance at least $2$ to
$\gamma_0$. Again, since $g$ is an eigenfunction to the eigenvalue
zero, $\gamma_2$ must have a neighbour $\gamma_3$ with $g(\gamma_3) >
0$. Therefore, $\gamma_3 \in \supp_+(g) \backslash \{ \gamma_0 \}$,
which is a contradiction.

So we proved $g(\gamma_0) = 0$. Let $\gamma' \in \supp(g)$. Since
$G$ is $3$-connected, there are three pairwise disjoint paths $P_1,
P_2, P_3$, connecting $\gamma_0$ with $\gamma'$. Without loss of
generality, we can assume that the path $P_i$ contains the vertex
$\gamma_0 s_i$. Starting in $\gamma_0 s_i$ and following the path
$P_i$ in direction $\gamma'$, let $\gamma_i \in P_i$ be the first
vertex with $g(\gamma_i) = 0$ and $\gamma_i$ adjacent to
$\supp(g)$. Since $g$ is an eigenfunction, $\gamma_i$ must be adjacent
to both $\supp_+(g)$ and $\supp_-(g)$. Now, contract $\supp_+(g)$ and
$\supp_-(g)$ to single vertices, denoted by $v_+$ and $v_-$ (which is
possible since both sets are connected, by Proposition
\ref{prop:courant}) and contract also the parts of the paths
$P_1,P_2,P_3$ from $\gamma_0$ to $\gamma_i$, and remove all other
vertices on which $g$ vanishes. The resulting graph is planar and
contains $K_{3,3}$ as a subgraph (where one set of vertices are
$\gamma_1,\gamma_2,\gamma_3$ and the other set are
$\gamma_0,v_+,v_-$), which is a contradiction. \EPf

Proposition \ref{prop:A3B3H3} follows now immediately from Proposition
\ref{prop:kassabov} and Corollary \ref{cor:mult3}.

\begin{rmk}
  Let $X \in {\rm int}(\Delta)$, and $\E$ be the eigenspace of $P_X$
  to the eigenvalue $\lambda_1(X)$. The above arguments show that, for
  all $\gamma \in \Gamma$, the maps $\Theta_\gamma: \E \to \R^3$
  (given by \eqref{eq:theta}) are {\em bijective}. This fact is
  equivalent to a particular transversality property of $P_X$, the
  so-called {\em Strong Arnold Hypothesis} (for the precise definition
  see, e.g., \cite{CdV91} or \cite{LSch-99}). The Strong Arnold
  Hypothesis played a crucial role in the proof that Colin de
  Verdi{\'e}re's graph invariant is monotone with respect to taking
  minors.
\end{rmk}

\section{Proofs of the results about the Archimedean solids} 
\label{sec:arch46x}

Before we present the proofs of Theorems \ref{thm:seceigvalarch} and
\ref{thm:threecurves}, let us mention that the full spectra of the
canonical Laplacians of the Archimedean solids were calculated in
\cite{ST-98}.

\smallskip

{\bf Proof of Theorem \ref{thm:seceigvalarch}:} 
Let $\Gamma \in \{A_3, B_3, H_3\}$ and $G=(V,E)$ be the Cayley graph
associated to $\Gamma$ with respect to the canonical set $S =
\{s_1,s_2,s_3\}$ of generators. Recall that $G$ is the one-skeleton of
the Archimedean solids $(4,6,6)$, $(4,6,8)$ and $(4,6,10)$,
respectively.

Then we have
\begin{equation} \label{eq:MM1}
M = \begin{pmatrix} 1 & 0 & -\frac{1}{2} \\ 0 & 1 & -\frac{\eta}{2} \\
  -\frac{1}{2} & -\frac{\eta}{2} & 1 \end{pmatrix}, \qquad M^{-1} =
\frac{1}{\rho} \begin{pmatrix} 1+\rho & \eta & 2 \\
  \eta & 3 & 2 \eta \\ 2 & 2 \eta & 4 \end{pmatrix} 
\end{equation}
and $V^2 = \frac{\rho}{4}$, where $\rho$ and $\eta$ are given as in
the following table:

\smallskip

\begin{center}
  \begin{tabular}{l|c|c|c} $\Gamma$ & $A_3$ & $B_3$ & $H_3$ \\ \hline
    $\eta$ & $1$ & $\sqrt{2}$ & $\varphi = \frac{1+\sqrt{5}}{2}$ \\
    $\rho=3-\eta^2$ & $2$ & $1$ & $2-\varphi$
\end{tabular}
\end{center}

\smallskip

Let $p = \alpha p_1 + \beta p_2 + \gamma p_3$ be a general
point in the spherical fundamental domain $\F_0 \subset S^2$. Choosing
$\mu=1$ and using \eqref{eq:calcX'} and \eqref{eq:MM1}, we obtain
\begin{equation} \label{eq:xyz'}
  \begin{pmatrix} x' \\ y' \\ z' \end{pmatrix} = \frac{1}{2\sqrt{\rho}} 
  \begin{pmatrix}
  1+\rho+\frac{\eta \beta +2 \gamma}{\alpha}\\ 3 +
  \frac{\alpha+2\gamma}{\beta} \eta \\ 4 + \frac{2\alpha+2\eta\beta}{\gamma} 
  \end{pmatrix} 
\end{equation}
and $\lambda'=x'+y'+z'-\sqrt{\rho}$. $\Psi_\Delta(p)$ and $\Psi_\lambda(p)$
are then given by the expressions in \eqref{eq:PsisP}.

Since the lengths of the Euclidean edges are given by $\Vert p -
\sigma_1(p) \Vert = 2 \alpha V$, $\Vert p - \sigma_2(p) \Vert = 2
\beta V$ and $\Vert p - \sigma_3(p) \Vert = 2 \gamma V$ (see Corollary
\ref{cor:coxgroup}), there is only one point $p_0 \in \F_0$ for which
all edges are of equal length, namely the choice $\alpha = \beta =
\gamma$. Using \eqref{eq:xyz'} in this case and calculating
$(x,y,z)=\Psi_\Delta(p_0)$ and $\lambda=\Psi_\lambda(p_0)$ with the help of
\eqref{eq:PsisP} yields
\begin{equation} \label{eq:xyzmin}
(x,y,z) = \frac{1}{12+\rho+6\eta}(3+\rho+\eta,3+3\eta,6+2\eta) \quad
\text{and}\ \lambda = \frac{12+6\eta-\rho}{12+6\eta+\rho}. 
\end{equation}
By Theorem \ref{thm:main} and Proposition \ref{prop:A3B3H3}, this is
the only critical point of $\lambda_1: {\rm int}(\Delta) \to
(0,1)$. By Proposition \ref{prop:lambda1boundary}, $\lambda_1$ has a
global minimum in ${\rm int}(\Delta)$, which must therefore agree with
\eqref{eq:xyzmin}.
  
From Corollary \ref{cor:coxgroup} and Proposition \ref{prop:A3B3H3} we
conclude that $\lambda_1: {\rm int}(\Delta) \to (0,1)$ is analytic, and
we know from the proof of Proposition \ref{prop:lambda1boundary} that
$\lambda_1$ is convex. Assume that $\lambda_1$ would not be {\em
  strictly convex}. Then there would exist three different collinear
points $X_1, X_2, X_3 \in {\rm int}(\Delta)$ with $\lambda_1(X_1) =
\lambda_1(X_2) = \lambda_1(X_3)$. Convexity of $\lambda_1$ would force
$\lambda_1$ to be constant on the line segment bounded by the two
extremal points of $X_1, X_2, X_3$. Analyticity of $\lambda_1$ would
imply that $\lambda_1$ is constant along the whole line in $\Delta$
containing these three points. But this would lead to $\lambda_1(X_1)
= \lambda_1(X_2) = \lambda_1(X_3) = 1$, a contradiction to
$\lambda_1<1$ on the interior of $\Delta$.

In the case $\Gamma=H_3$, i.e., $(\eta,\rho) = (\varphi,2-\varphi)$,
we obtain
$$ X_0=(x,y,z) = \frac{1}{14+5\varphi}(5,3+3\varphi,6+2\varphi) \quad
\text{and}\ \lambda = \frac{10+7\varphi}{14+5\varphi}. $$ 
The corresponding spectral representation $\Phi_{X_0}$ agrees, up to
the factor $\frac{|\Gamma|}{3}$, with the orbit map
$\Phi(\gamma)=\gamma p_0$, by Corollary \ref{cor:coxgroup}, and is
therefore faithful.

Analogously, one easily checks that the choices $(\eta,\rho) = (1,2)$
and $(\eta,\rho) = (\sqrt{2},1)$ lead to the explicit values for
$(x,y,z)$ and $\lambda$, given in Remark \ref{rmk:minima}. \EPf

\smallskip

{\bf Proof of Theorem \ref{thm:threecurves}:} We only discuss the
Archimedean solid $(4,6,10)$ (i.e., $\Gamma = H_3$), the other solids
are treated analogously. 

Note, by the construction of $p_1,p_2,p_3$ in \eqref{eq:Pj}, that the
orbit $\Gamma p_1$ gives the vertices of an icosahedron. Up to a
scalar factor, $p_2$ points to the centre of a face of this
icosahedron and $p_3$ to the midpoint of an edge of the icosahedron,
and the orbits $\Gamma p_2$ and $\Gamma p_3$ are the vertices of an
dodecahedron $(5,5,5)$ and of an icosidodecahedron $(3,5,3,5)$,
respectively.  Moreover, it is easy to see that there are positive
constants $0 < c_0 < C_0$ such that $\alpha p_1 + \beta p_2 + \gamma
p_3 \in \F_0$, $\alpha, \beta, \gamma \ge 0$, implies $c_0 \le \alpha
+ \beta + \gamma \le C_0$.

Let $X_n = (x_n,y_n,z_n) \in {\rm int}(\Delta)$ be a sequence
converging to $(x,0,z) \in \Delta$ with $x,z > 0$. Then there are
constants $c_1, c_2 > 0$ with $c_1 x_n < z_n < c_2 x_n$. Let
$\Psi_\Delta^{-1}(X_n) = q_n = \alpha_n p_1 + \beta_n p_2 + \gamma_n p_3 \in
\F_0$. Our aim is to show $\alpha_n, \gamma_n \to 0$. From
\eqref{eq:xyz'} and \eqref{eq:PsisP}, we deduce that
\begin{eqnarray*}
  x_n &=& \frac{1}{F(\alpha_n,\beta_n,\gamma_n)} 
  \left( 3-\varphi+\frac{\varphi \beta_n + 2 \gamma_n}{\alpha_n} \right), \\
  y_n &=& \frac{1}{\beta_n F(\alpha_n,\beta_n,\gamma_n)} 
  (3 \beta_n + (\alpha_n + 2 \gamma_n) \varphi), \\ 
  z_n &=& \frac{1}{F(\alpha_n,\beta_n,\gamma_n)}
  \left( 4+\frac{2 \alpha_n + 2 \varphi \beta_n}{\gamma_n} \right),
\end{eqnarray*}
with 
$$ F(\alpha,\beta,\gamma) = 10-\varphi + 
\frac{\varphi \beta+2 \gamma}{\alpha} +
\frac{\alpha+2\gamma}{\beta}\varphi + 
\frac{2\alpha + 2 \varphi \beta}{\gamma}. $$ 

Since $y_n \to 0$ and $c_0 \le 3 \beta_n + (\alpha_n + 2\gamma_n)
\varphi \le 4 C_0$, we must have $\beta_n F(\alpha_n,\beta_n,\gamma_n)
\to \infty$.  This necessarily implies $\alpha_n \gamma_n \to
0$. Assume $\alpha_n$ converges to zero on a subsequence, on which
$\gamma_n$ does not converge to zero. Then
$F(\alpha_n,\beta_n,\gamma_n) \to \infty$ implies that $z_n$ converges
to zero on a finer subsequence and, since $c_1 x_n < z_n$, $x_n$ must
also converge to zero on this finer subsequence, contradicting to $x_n
+ y_n + z_n = 1$. This shows that both $\alpha_n, \gamma_n \to 0$,
i.e., $q_n$ converges to a multiple of $p_2$. By Corollary
\ref{cor:coxgroup}, the corresponding spectral representations
converge, up to a scalar factor, to the orbit map $\Phi(\gamma) =
\gamma p_2$, and $\Gamma p_2$ are the vertices of a dodecahedron. This
proves the convergence behaviour as $X_n$ converges to an interior
point of the bottom edge of the simplex $\Delta$ in Figure
\ref{fig:archconv1}. The converge behaviour to interior points of the
other two edges of $\Delta$ is proved analogously.

The curve $C_2$ is characterised by the property
$\alpha=\gamma$. Using this fact and the relation \eqref{eq:xyz'}, and
substituting $t=\frac{\alpha}{\beta}$ we obtain
$$ C_2 = \left\{ \frac{1}{3\varphi t^2+(14-\varphi) t+3\varphi)}
  \begin{pmatrix} (5-\varphi)t+\varphi\\
    3\varphi t^2 + 3t \\ 6t+2\varphi \end{pmatrix} \mid t \in
  (0,\infty) \right\} \subset \Delta. $$ 
Note that $t \to \infty$ implies $\beta \to 0$, which means that the
Euclidean edges between the $4$-gons and the $10$-gons shrink to zero
and the corresponding spectral representations converge to equilateral
realisations of the buckeyball $(5,6,6)$ (see Figure
\ref{fig:transspecrep}). The convergence behaviour along the other
curves $C_1, C_3$ is proved analogously. \EPf

\bibliographystyle{plain}
\bibliography{biblio}

\end{document}